\documentclass[12pt]{article}

\usepackage[all]{xypic}
\usepackage{amsmath}
\usepackage{amsfonts}
\usepackage{amssymb}
\usepackage{amsthm}

\theoremstyle{plain}
\newtheorem{thm}{Theorem}[section]

\newtheorem{crl}[thm]{Corollary}

\newtheorem{defn}[thm]{Definition}

\newtheorem{lmm}[thm]{Lemma}

\theoremstyle{remark}
\newtheorem{exm}[thm]{Example}

\newtheorem{rmk}[thm]{Remark}

\newtheorem{note}[thm]{Note}

\newcommand {\mr}{\mathrm}

\newcommand {\mb}{\mathbb}
\newcommand {\Z}{\mb Z}
\newcommand {\R}{\mb R}
\newcommand {\C}{\mb C}

\newcommand {\N}{\mb N}
\newcommand {\tP}{{\mb P}_\tro}
\newcommand {\tro}{\mathrm{trop}}

\newcommand {\colim}{\textrm{colim}\ }

\begin{document}

\title{Tropicalisation for topologists}

\author{Hadi Zare}
\date{}

\maketitle

\abstract{We consider the problem of translating notions from
classical topology to tropical language. We consider the tropical
projective and Grassmannian spaces. We give a fairly easy
classification of the projective gadget whereas the Grassmannians
seems rather more difficult. We also consider the notion of the
tropical matrices, and define a variant of tropical orthogonal
matrices. We completely determine $Gl(\R_\geqslant,n)$ and
$O(n)^\tro$. We also give some results on the idempotent elements of
the tropical matrix algebra $M_{n,n}(\R_\geqslant)$. Such a notion
will be important in the related bundle theory. We note that
tropical phenomena have been studies by algebraic geometers, and our
work here may overlap with them and there is no claim on the
originality of these notes, neither no claim on them being
unoriginal.}

\tableofcontents

\section{Motivation, Basics}
The $min$-algebra, $\N^\infty=\N\sqcup\{+\infty\}$ together with the
binary operation of taking minimum is our main motivation.
\begin{defn}
By a tropical space, we mean a pointed set $(T,\infty)$ together
with a $minimum$-function $\min=\min_T:T\times T\to T$ satisfying
$$\min (t,\infty)=\min(\infty,t)=t\textrm{ for all }t\in T,$$
and $\mit(t,t)=t$ for all $t\in T$.
\end{defn}
Notice that under this operation, $T$ turns into a monoid, whose
each element is idempotent. The tropical space $T$ is called
commutative if $\min(t,t')=\min(t',t)$ for all $t,t'\in T$. Next, we
make it clear what we mean by a tropical module, tropical algebra,
etc. Let $R$ be a commutative ring. By a left tropical $R$-module
$T$ we mean a tropical set $T$ of which has the structure of a left
$R$-module, and the $R$-action fixes the $\infty$, i.e.
$r\infty=\infty$. Similarly, we can define a tropical algebra $T$ to
be a tropical space $T$ together with a multiplication which is
compatible with the tropical structure.\\

Notice that we may consider the category of these objects, say
category of tropical $R$-modules together with the forgetful functor
$$\textrm{RMod}_\tro\longrightarrow \textrm{RMod}.$$

Next, we show how it is possible to construct examples of these
spaces. This of course reduces to find reasonable ways of defining
the $min$-function.

\begin{note}
We note that it is possible to have a dual tropical structure given
by taking \textit{maximum} of two given elements. We will use this
in some of our important examples.
\end{note}

\subsection{Examples}

\subsubsection{Tropicalising normed spaces}
Let us write $\R_{\geqslant}$ for the half line of nonnegative real
numbers. Let $(M,|\ |)$ be a space with a norm function $|\
|:M\to\R_{\geqslant}$. Notice that any norm gives rise to a metric,
and any metric space $(M,d)$ with a chosen point $0\in M$ can be
turned into a normed space by setting $|v|=d(v,0)$ for any $v\in M$.
Having a norm define $m:M\to \R_{>}$ by
$$m(a,b)=\min(|a|,|b|)$$
where the right hand side is the minimum function of the real line.
Notice that if we choose $M=\R$ with the usual the metric on $\R$,
then $m(a,b)=\min(a,b)$. Now, we may define $\min=\min_M:M\times
M\to M$ by $\min(a,a)=a$ and
$$\min(a,b)=\left\{ \begin{array}{ll}
                    a&\textrm{if }m(a,b)=|a|,\\
                    b&\textrm{if }m(a,b)=|b|.
                    \end{array}\right.$$
For a given normed space $M$, we define \textit{tropicalisation} of
$M$ by $M^\tro=M\sqcup\{\infty\}$ with generalised $min$-function
satisfying
$$\min(m,\infty)=\min(\infty,m)\textrm{ for all }m\in M$$
and $\min(\infty,\infty)=\infty$. Hence we obtain, a tropicalisation
functor
$$Trop:\textrm{Normed-Space}\longrightarrow \textrm{Trop-Space}.$$
Notice that category of normed spaces contains the category of
finite dimensional vector spaces inside it.\\
The fact that the tropicalisation of $M$ depends on the norm, or say
the metric, shows that this can be used as an invariant of metric
spaces, and not necessarily as an invariant of topological spaces.
The reason is that it is possible to choose different metrics,
yielding the same topology on a space.
\begin{note}
We note that the tropicalisation functor in fact is defined from the
category of totally ordered sets into the category of tropical sets,
as the notion of order naturally tells how to choose smaller element
among two given ones.
\end{note}

\subsubsection{Tropicalising Modules with projections}
Notice that on the real line one has
$$\min(a,b)=\frac{a}{2}+\frac{b}{2}-\frac{1}{2}|a-b|.$$
We use this as analogy to define $\min$ on a class of modules.
Suppose we have a ring $R$ which includes $\frac{1}{2}$, say
$R=\Z[\frac{1}{2}], \mathbb{Q},\R,\C$, etc. Let $M$ be an $R$-module
together with a projection, i.e. a mapping $\tau:M\to M$ such that
$\tau^2=\tau$. We then define
$$\min{_M^\tau(a,b)}=\frac{1}{2}(a+b-\tau(a-b)).$$
This then defines a $min$-function on $M$. Note that one may
restrict to choose $\tau$ as an automorphism of $R$-modules. But for
the moment we shall not put this restriction. We can consider the
tropicalisation of $M^\tro_\tau=M\sqcup\{\infty\}$ with generalised
$min$-function defined as
$$\min{_{M^\tro_\tau}(m,\infty)}=\min{_{M^\tro_\tau}(\infty,m)}=m\textrm{ for all }m\in M,$$
and $\min{_{M^\tro_\tau}(\infty,\infty)}=\infty$. This then defines
the tropicalisation functor
$$Trop:\textrm{Rmod}^{\textrm{pro}}\longrightarrow \textrm{Trop-Space}$$
where the left hand side is the category of $R$-modules with
projections. If we give structure of a monoid to $M^\tro_\tau$ by
$\min{_{M^\tro_\tau}}$ and define the multiplication on
$M^\tro_\tau$ by generalising the addition operation of $M$, then we
obtain an example of a tropical algebra. More precisely, we define
$\oplus,\odot$ by
$$\begin{array}{ll}
m\oplus  m'=\min{_{M^\tro_\tau}(m,m')} \textrm{ for all }m,m'\in M^\tro_\tau,\\ \\
m\odot m'=\left\{ \begin{array}{ll}
                    m+m'   &\textrm{ if }m,m'\in M\\
                    \infty &\textrm{ if }m=\infty,\textrm{ or }m'=\infty.
                    \end{array}\right.
\end{array}$$
Notice that here $\odot$ is commutative. We usually will drop
$\odot$ from our notation. This then gives structure of a semi-ring
to $(M^\tro_\tau,\oplus,\otimes)$. Observe that the inclusion
$$M\to M^\tro_\tau$$
acts like the \textit{exponential map} as it sends addition to
multiplication. If we choose $\tau$ to be an automorphism, then we
can give structure of an $R$-module to to $M^\tro_\tau$ by requiring
that the action fixes $\infty$, and on the other points of
$M^\tro_\tau\setminus\{\infty\}$ acts in the same way as on $M$.
This then makes the inclusion map $M\to M^\tro_\tau$ into a map of $R$-modules.\\

Finally, we note that it is possible to define a ``dual'' tropical
structure on $M$. For motivation, notice that in $\R$, we have
$$\max(a,b)=\frac{1}{2}(a+b+|b-a|).$$
One then may fix $\tau:M\to M$ with $\tau^2=1$ and define
$$\max{_M^\tau(a,b)}=\frac{1}{2}(a+b+\tau(a-b)).$$
It is then possible to perfom as we did with the $min$-function.

\section{Tropical Euclidean and Grassmannian Spaces}
I like to approach these topics with hat of a topologist on. These
are quite familiar objects, and it will be interesting how they look
like geometrically.\\
In topology, there are different approaches to define Grassmannian
spaces. Let us to work with $\R$ for a moment. We consider $\R$ with
its usual addition and multiplication. The Grassmannian
$G_k(\R^{n+k})$ may be viewed as the space of $k$-dimensional
subspaces of $\R^{n+k}$ which can be described as quotient of
another space, namely the Steifel space $V_k(\R^{n+k})$ of
$k$-frames in $\R^{n+k}$ together the quotient map
$$V_k(\R^{n+k})\to G_k(\R^{n+k}).$$
Another way of defining $G_k(\R^{n+k})$ is to consider it as the
quotient
$$\frac{O(n+k)}{O(n)\times O(k)}$$
where $O(n)$ is the orthogonal group, the group of isometries of
$\R^n$.\\
The key element, in any of these approaches, is the type of
algebraic structure of $\R^n$ is a ``set with addition'' together
with its structure as an $\R$-module. On the other hand, the
standard $\R$-module structure on $\R^n$ is obtained by using the
multiplication operation $\R\times\R\to\R$.\\

We consider the tropicalisation of the real line and determine a
model for it, which is more familiar to a topologist. We then
consider the ways that one can define the multiplication operation
on the tropical line, have granted that the addition is given by the
tropical structure.

\subsection{A Model for $(\R^\tro)^n$}
To begin with, notice that $\R^\tro$ as a set is in a one to one
correspondence with $\R_{\geqslant}$ where by the latter we mean the
nonnegative real numbers. We fix such a correspondence as following.
First, notice that there is a homeomorphism of topological spaces
$f:\R_{>}\to\R$ defined by
$$f(x)=-\ln x$$
where $\R_{>}$ stands for the open half line of positive real
numbers. We extend this to $f:\R_{\geqslant}\to\R^\tro$ by setting
$$f(0)=\infty.$$
This then enables us to get a one to one correspondence
$$f^{\times n}=(f,\ldots,f):\R_{\geqslant}^n\to(\R^\tro)^n$$
where $X^n$ denotes $n$-fold Cartesian product of $X$ with itself.
Note that it is quite straightforward to see that the category of
tropical spaces is closed under the Cartesian product.
\begin{note}
For a topologist the object $\R_{\geqslant}^n$ is a familiar one.
This is the space that is used as the model space to define
\textit{manifolds with corners}, where
$$\partial\R_{\geqslant}^n=\R_{\geqslant}^n-\mathrm{interior}(\R_\geqslant^n)$$
corresponds to different types of (potential) singularities that
such a manifold can have, whereas the interior corresponds to the
smooth parts of such manifolds.
\end{note}

\begin{rmk}
Notice that the set $\R^\tro$ has two components, hence the set
$(\R^\tro)^n$ will have $2^n$ components. We may write this as
$$(\R^\tro)^n=\coprod_{0\leqslant k\leqslant n}(\R^k)^{\sqcup {n\choose k}},$$
where by $(\R^k)^{\sqcup {n\choose k}}$ we mean ${n\choose k}$-fold
$\coprod$-product of $\R^k$ with itself. We can make this more
precise. For instanse, in $(\R^\tro)^2$ we have one copy of $\R^2$,
two copies of $\R$ and one copy of $\R^0$ in the following way:
$\R^2$ corresponds to itself, one copy of $\R$ corresponds to
$\{\infty\}\times\R$, another copy of $\R$ corresponds to
$\R\times\{\infty\}$, and $\R^0$ corresponds to $(\infty,\infty)$.
Under the correspondence of $(\R^\tro)^2$ with $\R_{\geqslant}^2$ we
see that $\R^2$ corresponds to the interior of $\R_{\geqslant}^2$,
$\R\times\{\infty\}$ corresponds to the $x$-axis without the origin,
$\{\infty\}\times\R$ corresponds to the $y$-axis, and
$(\infty,\infty)$ corresponds to the origin.
\end{rmk}

Our next objective, is to define an action on $\R^\tro$ and make
sense of $k$-dimensional subspaces of $(\R^\tro)^{n+k}$. Here, we
have to choose which object we like to act on $\R^\tro$. We may
choose $\R$ to acts on $\R^\tro$, or we may choose $\R^\tro$ to act
on $\R^\tro$.
\subsubsection{A Tropical Structure on $\R^n_\geqslant$}
We consider the ``dual'' tropical structure on $\R_{\geqslant}$ by
using the taking the maximum of two given values. This operation has
a min element $0\in\R_{\geqslant}$ as its minimum. This then induces
a tropical structrue on $\R_{\geqslant}^n$ given
$$v\oplus v'=(\max(v_1,v'_1),\ldots,\max(v_n,v'_n))$$
where $v,v'\in \R_{\geqslant}^n$. In this case
$0\in\R_{\geqslant}^n$ correspond to the minmum element where its
image under $f^{\times n}:\R_{\geqslant}^n\to(\R^\tro)^n$
corresponds to maximum element, i.e. $\infty$. Notice that $g$ and
$f$ both are decreasing functions, hence they respect the tropical
structure. Observe that the positive real line $\R_>$ is a group
under multiplication. In fact $\R_\geqslant$ is the tropicalisation
of $\R_>$ viewed as a group under multiplication. This then enables
us define the multiplication
$\R_\geqslant\times\R_\geqslant\to\R_\geqslant$ to be the usual
multiplication of two real numbers, i.e.
$$(r,s)\longmapsto rs.$$
It is straightforward to see that $f$ and $g$ become maps of
tropical (semi)rings.\\

Our next goal is to make it clear what we mean by the action of real
numbers on these spaces.

\subsubsection{Flows on $\R_{\geqslant}^n$}
The one to one correspondence
$$f^{\times n}:\R_{\geqslant}^n\to(\R^\tro)^n$$
motivates us, and provides us with a tool, to define an action of
$\R_{\geqslant}$ on $\R_{\geqslant}^n$. Let us write $g$ for the
inverse of $f$, i.e. $g:\R^\tro\to\R_{\geqslant}$ is given by
$$\begin{array}{lll}
g(x)     &=&e^{-x}\\
g(\infty)&=&0.
\end{array}$$
Recall from previous section that we have tropicalisation of any
module with projection. In this case, $\R^\tro$ is the same as
tropicalisation of $\R$ with $\tau$ given by the norm function. This
then shows that it is possible to have a multiplication
$\odot:\R^\tro\times\R^\tro\to\R^\tro$ given by
$$t\odot t'=\left\{\begin{array}{ll}
                   t+t'   &\textrm{if }t,t'\in\R,\\
                   \infty &\textrm{if }t=\infty,\textrm{ or
                   }t'=\infty.
                   \end{array}\right.$$
We use this to define the action of the additive group $\R$ on
$\R^\tro$ as $\R\times\R^\tro\to\R^\tro$ given by
$$(r,t)\longmapsto r\odot t.$$
We use $g$ to define an action of
$\R\times\R_{\geqslant}\to\R_\geqslant$ by requiring that $g$
respects this action, i.e. $g$ has to be a map of $\R$-modules. This
induces the action $\R\times\R_\geqslant\to\R_\geqslant$ given by
$$(r,v)\longmapsto e^{-r}v.$$
Notice that here we consider the action of $(\R,+)$ on
$\R_\geqslant$. It is again clear that the mappings $f,g$ respect
these actions. Notice that it is then quite clear that how to define
the corresponding actions of $\R$ on $(\R^\tro)^n$ and
$\R_\geqslant^n$ respectively. This is done by defining the action
component-wise. We investigate this in the next part, where we look
at analogous of $k$-planes in these spaces.\\

\begin{note}
Finally we explain the tittle for this section. The word
\textit{flow}, for a differential topologist, reminds the action of
the real numbers under addition on a (smooth) manifold. The above is
only one particular flow that we use, and let call it the ``standard
flow''
on $\R_\geqslant$.\\
It is possible to consider a more general setting. Notice that we
have defined the mappings $f$ and $g$ in a way that they respect the
tropical (semi)ring structures. Hence, having any action of the
additive group $(\R,+)$ on $\R_\geqslant$ will determine a
corresponding action on $\R^\tro$, hence on $\R_\geqslant^n$ and
$(\R^\tro)^n$ respectively. Note that in the case of the standard
flow the point $0\in\R_\geqslant^n$ corresponds to a singular point
of the flow.
\end{note}

\subsection{Tropical Projective Spaces} The study of tropical
projective spaces, i.e. the space of lines in the tropical Euclidean
spaces, seems to be the first natural step in attempt to understand
the tropical Grassmannian spaces. It turns out that, using our model
for $(\R^\tro)^n$, it is an easy task to identify $\tP^n$ as a set
where $\tP^n$ is the set of all lines in $(\R^\tro)^{n+1}$ which
``pass'' through the origin. We state the result.
\begin{lmm}
There is a one to one correspondence
$$\tP^n\to\Delta^n$$
where $\Delta^n$ is the standard $n$-simplex, i.e.
$$\Delta^n=\{(x_1,\ldots,x_{n+1})\in\R^{n+1}:x_i\geqslant 0,\sum x_i=1\}.$$
\end{lmm}
This is quite straightforward to see, when we use
$\R^{n+1}_\geqslant$ as our model for $(\R^\tro)^{n+1}$. For
instance, let $n=1$. Then we need to identify all lines in
$\R^{n+1}_\geqslant$. Notice that in $\R^{n+1}_\geqslant$ any point
together with the ``origin'' determines a line. Let $(a,b)\in
\R^2_\geqslant$. Then the line passing through this point and
``reaching'' the origin is determined by the orbit of this point
under the action of real line, i.e. all points
$e^{-r}(a,b)=(e^{-r}a,e^{-r}b)$ where $r\in\R$. Notice that $e^{-r}$
is nothing but a positive real number. Hence, the orbit of $(a,b)$
is the line passing $(a,b)$ and the origin. However, notice that
this will never reach the origin, and the origin will be a limit
point for this line when $r$ tends to $+\infty$. This then shows
that the set $\{(a,b):a>0,b>0,a+b=1\}$ is in one to one
correspondence with these line. We need to identify the lines that
correspond to the cases with $a=0$ and $b=0$. The $x$-axis and the
$y$-axis then give these two end points. Notice that this look like two point
compactification of the real line.\\
For cases $n>2$ a similar approach gives the result. We only note
that the boundary of the simplex will correspond to the lines on
$\partial\R^{n+1}_\geqslant$ where as its interior points correspond
to $\mr{interior}(\R^{n+1}_\geqslant)$.
\begin{rmk}
The tropical projective space does not seem to be a tropical space,
however its correspondence an standard simplex may be an evidence
that it is a kind of a variety(?!).
\end{rmk}

\subsection{Subspaces in Tropical Euclidean Spaces}
We like to look at the tropical version of the Grassmannian spaces.
We note that there is some work on this from an algebric-geometric
point of view such as\\ \\
David Speyer, Bernd Sturmfels \textit{The tropical Grassmannian}
Adv. Geom. Vol.4 No.3 pp389--411, 2004\\

But I am not aware of the contents of this work, so I don't make any
comment. We choose to work with our model $\R^n_\geqslant$ to study
the Grassmannian objects. The Grassmannian space $G_k(\R^{n+k})$ is
the set of $k$-planes in $\R^{n+k}$. We take this approach and look
at the $k$-planes in $\R^{n+k}_\geqslant$.

\subsubsection{Linear independence in the tropical sense}
In the Euclidean space, it is quite straightforward to identify the
$k$-dimensional subspaces of $\R^{n+k}$, i.e. the space spanned by
$k$ linear independent vector $\{v_1,\ldots,v_k:v_i\in\R^{n+k}\}$.\\

In tropical space $(\R^\tro)^{n+k}$ two vectors
$(t_1,\ldots,t_n),(t'_1,\ldots,t'_n)$ are linearly dependent in the
tropical sense if and only if there exist $r\in\R$ such that
$$t_1+r=t'_1,\ldots,t_n+r=t'_n.$$

However, our model is much easier to use. More precisely, two
vectors $(t_1,\ldots,t_n),(t'_1,\ldots,t'_n)\in\R^{n+k}_\geqslant$
are linearly dependent if there is a real number $r$ such that
$$(t_1,\ldots,t_n)=e^{-r}(t'_1,\ldots,t'_n)$$
i.e.
$$t_1=e^{-r}t'_1,\ldots,t_n=e^{-r}t'_n.$$
This tells that two vectors in $\R^{n+k}_\geqslant$ are linear
independent if one is not multiple of the other one, similar to
the notion of the linear independent in the Euclidean sense. However, the notion
of linear combination is quite different in two cases.\\
Next, we have two identify what is meant by the space spanned by $k$
linear independent vectors $v_1,\ldots,v_k\in\R^{n+k}_\geqslant$ in
the tropical. The addition in $\R^{n+k}_\geqslant$ defined in
previous sections shows that spanning in tropical sense is given by
$$\mr{Span}^\tro\{v_1,\ldots,v_k\}=\mr{interior}(\mr{Cone}\{v_1,\ldots,v_k\})$$
where $\{v_1,\ldots,v_k\}$ is an arbitrary set of $k$ vectors in
$\R^{n+k}_\geqslant$ and the cone $\mr{Cone}\{v_1,\ldots,v_k\}$ is
the cone taken in usual Euclidean sense. For instance, for $k=2$ the
cone on two vectors is the area between the two lines determined by
two vectors. By a $k$-subspace $C\subset\R^{n+k}_\geqslant$ we mean
a cone which is span of $k$ independent vectors, i.e.
$$C=\mr{Span}^\tro\{v_1,\ldots,v_k\}$$
where $\{v_1,\ldots,v_k\}$ is a linearly independent set.\\
Notice that there is not a precise notion of dimension here. The
reason is that not every point in $\R^{n+k}_\geqslant$ is a linear
combination of finite number of vectors. The reason for this lies in
the way that the tropical addition, and scalar multiplication on
$\R^{n+k}_\geqslant$ are defined. The following observation provides
us with a framework to look at this.

\begin{lmm}
Let $u_1,\ldots,u_n\in\R^n$ denote the standard Euclidean basis
elements, i.e. $u_1=(1,0,\ldots,0),\ldots,u_n=(0,\ldots,0,1)$. We
then have
$$v\in\mr{interior}(\R^n_\geqslant)\Longleftrightarrow v\in\mr{Span}^\tro\{u_1,\ldots,u_n\},$$
where
$\mr{interior}(\R^n_\geqslant)=\R^n_\geqslant-\partial\R^n_\geqslant$.
The space $\partial\R^n_\geqslant$ is characterised by
$$(x_1,\ldots,x_n)\in \partial\R^n_\geqslant\Longleftrightarrow
x_t=0\textrm{ for some }1\leqslant t\leqslant n.$$
\end{lmm}

\begin{rmk}
We note that according to this lemma, the space $\R^n_\geqslant$ is
not finitely generated as an $(\R,+)$-set. The reason for this is
that we don't have $\infty$ in $\R$. Later on, we will consider the
action of $\R^\tro$ on this set, where $\R^n_\geqslant$ becomes an
$\R^\tro$-module.
\end{rmk}

Notice that a vector $v$ is in $\mr{interior}(\R^n_\geqslant)$ if
all of its component, written in the usual Euclidean basis, are
positive. Observe that in particular, if we choose any vector, $u_i$
such as in the above lemma, then under the correspondence
$\R^n_\geqslant\to(\R^\tro)^n$ we can see that
$$\mr{Span}^\tro\{u_i\}\simeq\R.$$
Moreover, let us write
$$\mr{Span}^\tro\{\widehat{u}_i\}\simeq\{\infty\}.$$
The above lemma together with the notation that just introduced our
allows us to formally rewrite the decomposition of $\R^n_\geqslant$
corresponding to the one given by Remark 2.2. The result reads as
following.
\begin{crl}
The space $\R^n_\geqslant$ has the following decomposition as
$$\begin{array}{lll}
\mr{Span}^\tro\{u_1,\ldots,u_n\}\sqcup\\ \\
\mr{Span}^\tro\{u_1,\ldots,u_{n-1},\widehat{u}_n\}\sqcup\cdots\sqcup \mr{Span}^\tro\{\widehat{u}_1,u_2,\ldots,u_n\}\sqcup\\ \\
\mr{Span}^\tro\{u_1,\ldots,u_{n-2},\widehat{u}_{n-1},\widehat{u}_n\}\sqcup\cdots\sqcup\mr{Span}^\tro\{\widehat{u}_1,\widehat{u}_2,u_3,\ldots,u_n\}\sqcup\\
\\
\sqcup \cdots \sqcup\\ \\
\mr{Span}^\tro\{u_1,\widehat{u}_2,\ldots,\widehat{u}_n\}\sqcup\cdots\sqcup\mr{Span}^\tro\{\widehat{u}_1,\widehat{u}_2,
\ldots,\widehat{u}_{n-1},u_n\}\\ \\
\sqcup \{(0,0,\ldots,0)\}.
               \end{array}$$
Here $\widehat{u}_i$ means that the vector $u_i$ is not in the set.
Moreover, under the correspondence $$\R^n_\geqslant\to(\R^\tro)^n$$
the space
$\mr{Span}^\tro\{u_1,\ldots,u_{i-1},\widehat{u}_i,\ldots,\widehat{u}_j,u_{j+1}\ldots
u_n\}$ maps to
$$\begin{array}{l}
\mr{Span}^\tro\{u_1\}\times\cdots\times\mr{Span}^\tro\{u_{i-1}\}\times\mr{Span}^\tro\{\widehat{u}_i\}\cdots\times\\ \\
\mr{Span}^\tro\{\widehat{u}_j\}\times\mr{Span}^\tro\{u_{j+1}\}\times\cdots\times\mr{Span}^\tro\{u_n\}\end{array}$$
which is the same as
$$\R\times\cdots\times\R\times\{\infty\}\times\cdots\times\{\infty\}\times\R\times\cdots\times\R$$
where for each $u_k$ in the spanning set we obtain a copy of $\R$ at
$k$th position, and for each $\widehat{u}_i$ we obtain a copy of
$\{\infty\}$ at the $i$th position.
\end{crl}

Finally, we have a little observation which will be important later
on when we consider the generalised tropical Grassmannians.
\begin{lmm}
Let $v_1,\ldots,v_k\in\R^n_\geqslant$ be linearly independent with
$k>1$. Let $v\in\mr{Span}^\tro\{v_{\alpha_1},\ldots,v_{\alpha_j}\}$
where $1\leqslant\alpha_1,\ldots,\alpha_j\leqslant k$ and $j<k$.
Then $v\not\in\mr{Span}^\tro\{v_1,\ldots,v_k\}$. In particular,
$v_i\not\in \mr{Span}^\tro\{v_{\alpha_1},\ldots,v_{\alpha_j}\}$.
\end{lmm}
This is easy to see, as if
$v\in\mr{Span}^\tro\{v_{\alpha_1},\ldots,v_{\alpha_j}\}$ then $v$
has to belong to boundary of the $k$-subspace determined by
$v_1,\ldots,v_k$. The result then follows from Corollary 2.7.

\subsubsection{Relating $G_k(\R^n_\geqslant)$ to configuration spaces}
We like relate the Grassmannian space $G_k(\R^n_\geqslant)$ to some
configuration spaces. The mapping fails to be an isomorphism. But it
at least provides a tool which presumably will help to analyse, and
calculate more sophisticated algebraic invariants of these spaces.

We start by recalling the analogous construction in homotopy. Let us
fix an arbitrary basis for $\R^n$. Let $V\subset\R^n$ be a
$k$-dimensional subspace. We then can choose a basis for it, say
$\{v_1,\ldots,v_k\}$. The fact that there are linearly independent
means that they give rise to $k$ distinct lines, hence defines a
mapping
$$G_k(\R^n)\to F(P^{n-1},k),$$
where given any set $X$ we define the set of configuration of $n$
point in $X$ as
$$F(X,n)=\{(x_1,\ldots,x_n):x_i\in X, i\neq j\Longrightarrow x_i\neq x_j\}.$$
The above mapping fails to be a homeomorphism as it is possible to
choose $k$ distinct points in $P^{n-1}$, or say $k$ distinct
vectors, which are not necessarily linearly independent.\\

We write $G_k(\R^n_\geqslant)$ for the set of all $k$-subspaces of
$\R^n_\geqslant$. Assume that we have a $k$-subspace
$C=\mr{Span}^\tro\{v_1,\ldots,v_k\}\subset\R^n_\geqslant$. Let
$l_1,\ldots, l_k$ be $k$ distinct lines determined by
$v_1,\ldots,v_k$ respectively. This determines a mapping
$$G_k(\R^n_\geqslant)\longrightarrow F(\tP^{n-1},k).$$
Notice that $\tP^{n-1}$ is the same as $\Delta^{n-1}$. The fact that
$v_1,\ldots,v_k$ are linearly independent in the tropical sense,
implies that non of these vectors falls into the cone generated by
the other ones. If we use $v_1,\ldots,v_k$ to denote those $k$
points in $\Delta^{n-1}$ that correspond to these lines, we obtain a
convex set, where here convex means convex as a subset of $\R^n$
with its usual metric.\\
On the other hand, if we choose any convex set in $\Delta^{n-1}$
with $k$ vertices we obtain $k$ vector in $\R^n_\geqslant$ which are
independent in the tropical sense. This completes the proof of the
following observation.
\begin{thm}
There is an isomorphism of sets
$$G_k(\R^n_\geqslant)\to F(\Delta^{n-1},k)^\mr{convex}$$
where $F(\Delta^{n-1},k)^\mr{convex}$ refers to a subset of
$F(\Delta^{n-1},k)$ whose points are in one to one correspondence
with convex subset of $\Delta^{n-1}$ with $k$ vertices.
\end{thm}

\subsection{$\R^n_\geqslant$ as an $\R^\tro$-module}
Recall from previous sections that $\R\times
\R^n_\geqslant\to\R^n_\geqslant$ gives $\R^n_\geqslant$ structure of
an $(\R,+)$-set. This action is not compatible when we consider the
field of real lines, with its usual addition and multiplication.
However, it is possible to obtain structure of an $\R^\tro$-module
on $\R^n_\geqslant$.\\
First, define $\R^\tro\times\R_\geqslant\to\R_\geqslant$ by
$$\begin{array}{ll}
(r,t) &\longmapsto e^{-r}t,\\
(+\infty,t)&\longmapsto 0.
\end{array}$$
Recall that $\R^\tro$ has a (semi)ring structure when regarded as
$(\R\cup\{+\infty\},\min,+)$ whereas $\R_\geqslant$ has the tropical
structure when regarded as $(\R_\geqslant,\max,\cdot)$ where $\cdot$
denotes the usual product. We then define the action
$\R^\tro\times\R^n_\geqslant\to \R^n_\geqslant$ to be the
component-wise action, i.e.
$$\begin{array}{ll}
(r,(t_1,\ldots,t_n))      &\longmapsto    (e^{-r}t_1,\ldots,e^{-r}t_n),\\
(+\infty,(t_1,\ldots,t_n))&\longmapsto(0,\ldots,0).
\end{array}$$
This definition is very similar to the previous one. It does not
change the notion of the linear independence. However, there is
slight difference in the notion of span. We write
$\mr{Span}^{\mr{Trop}}$ to distinguish it from $\mr{Span}^\tro$.
\begin{lmm}
Suppose $v_1,\ldots,v_k\in\R^n_\geqslant$. Then
$$\mr{Span}^{\mr{Trop}}\{v_1,\ldots,v_k\}=\mr{Cone}\{v_1,\ldots,v_k\}.$$
\end{lmm}
Hence, a slight change in the ground set acting on $\R^n_\geqslant$,
i.e. replacing $\R$ with $\R^\tro$ has the effect that it adds the
\textit{limit points} of a cone to it. As a corollary
$\R^n_\geqslant$ becomes finitely generated over $\R^\tro$. We have
the following.
\begin{crl}
Suppose $u_1,\ldots,u_n$ denote the usual basis elements for $\R^n$.
We then have
$$\R^n_\geqslant=\mr{Span}^{\mr{Trop}}\{u_1,\ldots,u_n\}.$$
Here, any point on $\partial\R^n_\geqslant$ will have tropical
coordinates which are formed of real numbers, and $+\infty$.
Moreover, $\{u_1,\ldots,u_n\}$ is the only set of vectors satisfying
this property, i.e. if there is any set of vectors
$\{v_1,\ldots,v_n\}$ such that
$$\R^n_\geqslant=\mr{Span}^{\mr{Trop}}\{v_1,\ldots,v_n\},$$
then each $v_i$ will be a re-scaling of of $u_j$ for unique
$1\leqslant j\leqslant n$.
\end{crl}
As an example, consider $\R^2_\geqslant$ and consider the point
$(1,0)$ which is on its boundary. The Corollary 2.7 implies that it
can not be written as any linear combination of two vectors in
$\R^2_\geqslant$ when regarded as an $\R$-set. However, as an
$\R^\tro$-module, we may write
$$(1,0)=e^0(1,0)+\infty(0,1)$$
i.e. as a vector in $\R^2_\geqslant$ the point $(1,0)$ may be
written as the column vector
$$\left [\begin{array}{c}
1\\
+\infty
\end{array}\right ].$$
\subsubsection{The space $G_k^\mr{Trop}(\R^n_\geqslant)$}
Likewise the space $G_k(\R^n_\geqslant)$ we define
$G_k^\mr{Trop}(\R^n_\geqslant)$ to be the set of all $k$-subspaces
in $\R^n_\geqslant$ when regarded as an $\R^\tro$-module. We say
$C\subseteq\R^n_\geqslant$ is a $k$-subspace if there are $k$
linearly independent vector $v_1,\ldots,v_k\in \R^n_\geqslant$ such
that
$$C=\mr{Span}^{\mr{Trop}}\{v_1,\ldots,v_k\}.$$
We may refer to $G_k^\mr{Trop}(\R^n_\geqslant)$ as the generalised
tropical Grassmannian space. Notice that there is a one to one
correspondence
$$G_k(\R^n_\geqslant)\longrightarrow G_k^\mr{Trop}(\R^n_\geqslant)$$
given by
$$C\longmapsto \overline{C}$$
where $\overline{C}$ denotes the closure of $C$, i.e.
$\overline{C}=C\cup\partial C$. The inverse mapping
$$G_k^\mr{Trop}(\R^n_\geqslant)\longrightarrow G_k(\R^n_\geqslant)$$
given by
$$C\longmapsto\mr{interior}{C}=C-\partial C.$$
Accordingly we obtain the following description of
$G_k^\mr{Trop}(\R^n_\geqslant)$.
\begin{thm}
There is a one to one correspondence
$$G_k^\mr{Trop}(\R^n_\geqslant)\longrightarrow F(\Delta^{n-1},k)^\mr{convex}.$$
\end{thm}

\begin{rmk}
Before proceeding further, we like to draw the reader's attention to
an essential difference between $G_k(\R^n_\geqslant)$ and
$G_k^\mr{Trop}(\R^n_\geqslant)$ in one hand and their Euclidean
analogous $G_k(\R^n)$ on the other hand. In the Euclidean space
$\R^n$ any set of $n$-linearly independent set will span $\R^n$,
however in $\R^n_\geqslant$ either as a $(\R,+)$ or as an
$\R^\tro$-module the only option for such a choice is provided by
the standard basis. Although, according to Lemma 2.6 as an
$(\R,+)$-set this it is not possible to generate all of
$\R^n_\geqslant$ in tropical sense.\\
For instance, consider $\R^3_\geqslant$ and let $C\in
G_k(\R^3_\geqslant)$ be defined as
$$C=\mr{Span}^\tro\{(1,0,0),(0,1,0),(1,1,2)\}.$$
This is a $3$-subspace in $\R^3_\geqslant$ and yet it is not equal
to $\R^3_\geqslant$. We note that all of those three vectors
defining $C$ are linearly independent. We can also consider to
$\overline{C}\in G_3^\mr{Trop}(\R^3_\geqslant)$ where
$C\neq\R^3_\geqslant$.\\
A consequence of this is that we may choose another vector
$v\in\R^3_\geqslant$ which does not belong to $C$, i.e the set
$$\{(1,0,0),(0,1,0),(1,1,2),v\}$$
is a linearly independent set in the tropical sense. This determines
a cone formed by $4$ vectors which can not be generated by any $3$
vectors. Hence we obtain a $4$-subspace in $\R^3_\geqslant$ giving
rise to an element of $G_4(\R^3_\geqslant)$.\\
In general, we may choose $k>n$ when we consider
$G_k(\R^n_\geqslant)$. In fact $k$ can be any arbitrary number.\\
\end{rmk}

It is quite interesting to see how a $k$-subspace in
$\R^n_\geqslant$ maps under the tropical isomorphism
$$\R^n_\geqslant\to(\R^\tro)^n.$$

Recall from previous sections that $(\R^\tro)^n$ is disjoint union
of $2^n$ copies of $\R^m$ with $0\leqslant m\leqslant n$. Now assume
that $C\in G_k(\R^n_\geqslant)$, i.e.
$$C=\mr{interior}(\mr{Cone}\{v_1,\ldots,v_k\})$$
where $v_1,\ldots,v_k\in\R^n_\geqslant$ are linearly independent(in
the tropical sense). Notice that in this case $C$ is given by the
interior of the cone, which implies that
$C\subset\mr{interior}(\R^n_\geqslant)$. We now that
$\mr{interior}(\R^n_\geqslant)$ maps to $\R^n\subset(\R^\tro)^n$.
Moreover, note that the mapping $f:\R^n_\geqslant\to(\R^\tro)^n$ is
continuous when restricted to $\mr{interior}(\R^n_\geqslant)$. This
implies that $C$ also maps into $\R^n\subset(\R^\tro)^n$.\\

Next, we consider $C\in G_k^\mr{Trop}(\R^n_\geqslant)$ and its image
under the $f:\R^n_\geqslant\to (\R^\tro)^n$. Notice that if $C\in
G_k^\mr{Trop}(\R^n_\geqslant)$ then it is a closed cone, where by
being closed we mean closed as a set in $\R^n_\geqslant$ viewed as a
topological space with its topology inherited from the standard topology on $\R^n$.\\
If $C\subset\mr{interior}(\R^n_\geqslant)$ then according to the
previous case, all of $C$ maps to only one component of $(\R^\tro)^n$, namely to $\R^n$.\\
Another possibility is that $C\cap\partial\R^n_\geqslant\neq\phi$.
In this case then the image of $C$ under
$f:\R^n_\geqslant\to(\R^\tro)^n$ will land in more than one factor
of $(\R^\tro)^n$ viewed as a disjoint union. The following provides
us with an example.
\begin{exm}
Consider the space $\R^3_\geqslant$, together with vectors
$v_1=(1,0,1)$ and $v_2=(0,1,1)$. This determines
$$C=\mr{Span}^\mr{Trop}\{v_1,v_2\}$$
as an element of $G_2^\mr{Trop}(\R^3_\geqslant)$. Let
$$L_1=\{(t,0,t):t>0\},\textrm{ }L_2=\{(0,t,t):t>0\},$$
i.e. $L_i\cup\{(0,0,0)\}$ is the line determined by $v_i$. It is
then clear that
$$\partial C=L_1\cup L_2.$$
Under the correspondence $f:\R^3_\geqslant\to(\R^\tro)^3$ we have
$$\begin{array}{lll}
f(L_1)&\subset&\R\times\{\infty\}\times\R\\
f(L_2)&\subset&\{\infty\}\times\R\times\R\\
f(0,0,0) & = & (\infty,\infty,\infty).
\end{array}$$
The image of $C$ under $f$ is an example of a
$\mathbf{2}$-subspace in $(\R^\tro)^3$.\\
\end{exm}

In general, suppose $C\in G_k^\mr{Trop}(\R^n_\geqslant)$ with
$C=\mr{Span}^\mr{Trop}\{v_1,\ldots,v_k\}$. If
$C\cap\partial\R^n_\geqslant\neq\phi$ then there are
$\alpha_i\in\{1,\ldots k\}$ with
$v_{\alpha_i}\in\partial\R^n_\geqslant$. Recall from Corollary 2.8
that $\R^n_\geqslant$ has a decomposition as
$$\begin{array}{lll}
\mr{Span}^\tro\{u_1,\ldots,u_n\}\sqcup\\ \\
\mr{Span}^\tro\{u_1,\ldots,u_{n-1},\widehat{u}_n\}\sqcup\cdots\sqcup \mr{Span}^\tro\{\widehat{u}_1,u_2,\ldots,u_n\}\sqcup\\ \\
\mr{Span}^\tro\{u_1,\ldots,u_{n-2},\widehat{u}_{n-1},\widehat{u}_n\}\sqcup\cdots\sqcup\mr{Span}^\tro\{\widehat{u}_1,\widehat{u}_2,u_3,\ldots,u_n\}\sqcup\\
\\
\sqcup \cdots \sqcup\\ \\
\mr{Span}^\tro\{u_1,\widehat{u}_2,\ldots,\widehat{u}_n\}\sqcup\cdots\sqcup\mr{Span}^\tro\{\widehat{u}_1,\widehat{u}_2,
\ldots,\widehat{u}_{n-1},u_n\}\\ \\
\sqcup \{(0,0,\ldots,0)\}.
               \end{array}$$
In this decomposition the first factor, i.e.
$\mr{Span}^\tro\{u_1,\ldots,u_n\}$ corresponds to
$\mr{interior}(\R^n_\geqslant)$, and the other factors correspond to
$\partial\R^n_\geqslant$. Hence, each of $v_{\alpha_i}$ will belong
to one, and only one, of the factors in the above decomposition. For
instance, assume
$$v_{\beta_1},v_{\beta_2},\ldots,v_{\beta_j}\in\mr{Span}^\tro\{u_1,\ldots,\widehat{u}_i,\ldots,\widehat{u}_j,\ldots,u_n\}=:S.$$
Notice that apart from $\{(0,0,\ldots,0)\}$, any other factor in the
above decomposition is an open set, when viewed as a subspace of
$\R^n$ with its usual topology. This then implies that
$$\mr{Span}^\mr{Trop}\{v_{\beta_1},v_{\beta_2},\ldots,v_{\beta_j}\}\subset\mr{Span}^\tro\{u_1,\ldots,\widehat{u}_i,\ldots,\widehat{u}_j,\ldots,u_n\}.$$
Applying Corollary 2.8 shows that
$\mr{Span}^\mr{Trop}\{v_{\beta_1},v_{\beta_2},\ldots,v_{\beta_j}\}$
maps into the factor of $(\R^\tro)^n$ corresponding to $S$. This
then completely determines where $C\cap\partial\R^n_\geqslant$ maps
under the tropical isomorphism $f:\R^n_\geqslant\to(\R^\tro)^n$.
This concludes our notes on the generalised Grassmannian tropical
spaces.

\subsubsection{$\mathbf{k}$-subspaces of $(\R^\tro)^n$.}
We like to analyse the those subspaces of $(\R^\tro)^n$ which are in
the images of $G_k^\mr{Trop}(\R^n_\geqslant)$ and map to more than
one factor in the disjoint union decomposition for $(\R^\tro)^n$.\\
We define what is meant by a $\mathbf{k}$-subspace in $(\R^\tro)^n$.
Recall that an example of a $\mathbf{2}$-subspace was given the
previous section.\\

In order to proceed, we need to fix an order on the disjoint union
decomposition for $(\R^\tro)^n$. Let $I=(i_1,\ldots,i_r)$ with each
of its entries belonging to $\{0,1\}$. By the set $(\R^\tro)^n_I$ we
mean a product of copies of the real line $\R$, and the singleton
$\{\infty\}$ in the following way: if $i_j=1$ then we have a copy of
$\R$, and if $i_j=0$ then we have a copy of $\{\infty\}$. For
example, for $n=2$, we have
$(\R^\tro)^2_{(1,0)}=\R\times\{\infty\}$. It is then clear that
$$(\R^\tro)^n=\bigsqcup_{I\in\{0,1\}^n}(\R^\tro)^n_I.$$
Moreover, let $|I|=\sum 2^{i_j}+1$. We then will refer to
$(\R^\tro)^n_I$ as the $|I|$-th factor of $(\R^\tro)^n$. This also
induces an order on these sequences (really the binary expansion of
positive natural numbers) by
$$I>J\Longleftrightarrow |I|>|J|.$$
This is the same as the the lexicographic order on the sequences $I$
and induces an order on the factors of $(\R^\tro)^n$ as following
$$(\R^\tro)^n_I\leqslant (\R^\tro)^n_J \Longleftrightarrow I\leqslant J.$$
Next, let $k>0$ and let $\mathbf{k}=(k_1,\ldots,k_{2^n})\in\Z^{2^n}$
be a sequence of nonnegative integers, with the most left nonzero
entry equal to $k$. Consider a collection of spaces
$\mathbf{K}=\{K_i:1\leqslant i\leqslant 2^n\}$ where $K_i$ is a
subspace of the $i$-th factor of $(\R^\tro)^n$ with
$$K_i=\mr{Span}^\tro\{v^i_1,\ldots,v^i_{k_i}\}$$
and $\{v^i_1,\ldots,v^i_{k_i}\}$ being a linear independent set in
the tropical sense. Moreover, we set the span of the empty set to be
the empty set. We call $\mathbf{K}$ as $\mathbf{k}$-subspace, if
there is a $k$-subspace $C\in G_k^\mr{Trop}(\R^n_\geqslant)$ such
that $f(C)$ maps into $(\R^\tro)^n$ with its image in different
factors of $(\R^\tro)^n$ being given by $K_i$'s.

\begin{rmk}
It is possible to give a more explicit account of the above
calculations. In order to do this, we need to label different
components of $\R^n_\geqslant$. It is similar to what we did in
above. Let $I$ be a sequence of length $n$ whose entries are either
$1$ or $0$. Let $u_1,\ldots,u_n$ denote the usual Euclidean basis
for $\R^n$. Let $(\R^n_\geqslant)_I\subset\R^n_\geqslant$ be given
by a product of copies of the open real half line $\R_>$ and the
singleton $\{0\}$, where for $i_j=1$ we have a copy of $\R_>$ at
$j$th position, and for $i_j=0$ we have have a copy of $\{0\}$ at
the $j$th position. For example, in the case of $n=2$ we have
$(\R^2_\geqslant)_{(1,0)}=\R_>\times\{0\}$ which is the $x$-axis
without $\{(0,0)\}$. We refer to $(\R^n_\geqslant)_I$ as the
$|I|$-th component of $\R^n_\geqslant$. Notice that
$(\R^n_\geqslant)_I=\mr{Span}^\tro\{e_1,\ldots,e_n\}$ where
$e_j=u_j$ if $i_j=1$ and $e_j=\widehat{u}_j$ if $i_j=0$. Hence,
according to Corollary 2.8 we have
$$\R^n_\geqslant=\bigcup_{I\in\{0,1\}^n}(\R^n_\geqslant)_I.$$
It is now evident that the $|I|$-th component of $(\R^n_\geqslant)$
maps homeomorphically onto the $|I|$-th factor of $(\R^\tro)^n$.\\

Now assume that $C\in G_k^\mr{Trop}(\R^n_\geqslant)$ with
$C\cap\partial\R^n_\geqslant\neq\phi$. Let
$C_I=C\cap(\R^n_\geqslant)_I$ be the $|I|$-th face of $C$. Then
$C_I$ maps homeomorphically into the $|I|$-th factor of
$(\R^\tro)^n$ under the tropical isomorphism
$$\R^n_\geqslant\to (\R^\tro)^n.$$
\end{rmk}

\subsection{Tropical Matrices, Tropical Orthogonal Matrices}
An alternative description of Grassmannian spaces is provided by
viewing them as the orbit space of two orthogonal groups acting on
another one. This is our main motivation study ``tropical matrices''.\\

Let $M_{m,n}(\R_\geqslant)$ be the set of all $m\times n$ matrix
whose entries are elements in $\R_\geqslant$. This set admits
structure of a semi-ring induces by the addition and multiplication
in $\R_\geqslant$. More precisely, for $A\in M_{m,n}(\R_\geqslant)$
let $A_{ij}$ denote the $(i,j)$ entry of $A$. Then for $A,B\in
M_{m,n}(\R_\geqslant)$ we define the addition, denoted with
$\oplus$, by
$$\begin{array}{lllll}
(A\oplus B)_{ij} & = & A_{ij}\oplus B_{ij} &=& \max(A_{ij},B_{ij}).
\end{array}$$
If $A\in M_{m,p}(\R_\geqslant)$ and $B\in M_{p,n}(\R_\geqslant)$,
then the multiplication, denotes with $\odot$, is defined by
$$\begin{array}{lllll}
(A\odot B)_{ij}  & = & \oplus_{k=1}^p A_{ik}\odot B_{kj} & =
&\max(A_{ik}B_{kj}:1\leqslant k\leqslant p).
\end{array}$$
Under these operations, the identity element for $\oplus$ is the
zero matrix, where as the identity matrix for $\odot$ is the
identity matrix. Moreover, notice that
$M_{m,n}(\R_\geqslant)$ is an $\R^\tro$-module.\\
One may ask whether if any matrix in $M_{n,n}(\R_\geqslant)$ is
invertible. The answer is positive, and an example is provided by
the set of all diagonal matrices whose diagonal does not have any
nonzero element. This is a ``universal example'' of such matrices as
the following theorem confirms, by giving a complete classification
of all such matrices.
\begin{thm}
Let $A\in M_{n,n}(\R_\geqslant)$ be a matrix with right
$\odot$-inverse with , i.e. there exists $B\in
M_{n,n}(\R_\geqslant)$ such that $A\odot B=I$. Then there exists
$\sigma\in\Sigma_n$ and a diagonal matrix $D=(D_1,\ldots,D_n)\in
M_{n,n}(\R_\geqslant)$ whose diagonal entries are nonzero, and
$$A=(D_{\sigma^{-1}(1)},\ldots,D_{\sigma^{-1}(n)}).$$
Here $\Sigma_n$ is the permutation group on $n$ letters. Moreover,
the matrix $B$ is determined by
$$B_{ij}=\left\{\begin{array}{ll}
                0&\textrm{if }A_{ji}=0\\
                \frac{1}{A_{ji}}&\textrm{if }A_{ji}\neq 0.
                \end{array}
                \right.$$
\end{thm}
Notice that according to this theorem any matrix with right inverse
also has a left inverse and they are the same. This is again
straight away to see this once we observe that
$$\begin{array}{lllll}
(A\odot B)_{ii}&=&\max(A_{ik}B_{ki}:1\geqslant k\geqslant n)&=&1,\\ \\
(A\odot B)_{ij}&=&\max(A_{ik}B_{kj}:1\geqslant k\geqslant
n)&=&0\textrm { \ for }i\neq j.
\end{array}$$
For instance, let $i=1$. The fact that $(A\odot B)_{11}=1$ implies
that there exists $k$ such that $A_{1k}=\frac{1}{B_{k1}}\neq 0$.
Combining $A_{1k}\neq 0$ together with $(A\odot B)_{ij}=0$ for
$i\neq j$ implies that $B_{kj}=0$ for $j\neq i$. Applying this for
each row of $A\odot B$ completes the proof.\\

Next, we consider the problem of determining tropical orthogonal
matrices. Recall that in the Euclidean case, the orthogonal matrix
$O(n)$ is an $n\times n$ matrix whose each column viewed has a
vector in $\R^n$ has unit norm, i.e. lies on the $(n-1)$-sphere
$S^{n-1}=\{x\in\R^n:||x||=1\}$, and is perpendicular to all other columns.\\
We may define the tropical inner product $\langle -,-
\rangle_\tro:\R^n_\geqslant\times\R^n_\geqslant\to\R^\tro$ by
$$
\langle v,w \rangle_\tro  =\oplus_{i=1}^nv_i\odot
w_i=\max(v_iw_i:1\geqslant i\geqslant n)
$$
with $v=(v_1,\ldots,v_n)$ and $w=(w_1,\ldots,w_n)$. In particular,
we have the tropical norm on $\R^n_\geqslant$ given by
$$||x||_\tro=(\max(x_i^2:1\leqslant i\leqslant n))^{1/2}.$$
In particular, the tropical circle is given by
$$\begin{array}{lll}
S^1_\tro & = & \{x\in\R^2_\geqslant:||x||_\tro=1\}\\
         & = & \{(1,x_2)\in\R^2:x_2\leqslant 1\}\cup\{(x_1,1):x_1\leqslant 1\}.
         \end{array}$$
More generally, the tropical $n$-sphere is give by
$$\begin{array}{lll}
S^n_\tro & = & \{(x_1,\ldots,x_{n+1})\in\R^{n+1}_\geqslant:||x||_\tro=1\}\\ \\
         & = & \bigcup_{i=1}^{n+1}\{(x_1,\ldots,x_{n+1})\in\R^{n+1}_\geqslant:x_i=1, j\neq i\Longrightarrow x_j\leqslant 1\}
         \end{array}$$
Moreover, these spaces admit a tropical structure.
\begin{lmm}
The tropical sphere $S^n_\tro$ together with the maximum operation,
inherited from $\R^{n+1}_\geqslant$, is a tropical space with the
identity element for this operation given by
$$(1,1,\ldots,1).$$
\end{lmm}

It is quite tempting to see what the analogous of the orthogonal
group will be. It is quite easy to determine form of such matrices.
The reason is provided with the following lemma.
\begin{lmm}
Let $A=(A_1,\ldots,A_n)\in M_{n,n}(\R_\geqslant)$ where $A_i$
denotes the $i$-th column of $A$. Suppose $||A_i||_\tro=1$ and
$\langle A_i,A_j\rangle_\tro=0$ for $i\neq j$. Then $A$ has the same
columns as the identity matrix.
\end{lmm}
Let $O(n)^\tro$ be the set of all matrices identifies by the above
lemma, i.e. set of all tropical orthogonal matrices.
\begin{lmm}
Let $\Sigma_n$ denote the permutation group on $n$ letters. Let the
action $\Sigma_n\times M_{m,n}(\R_\geqslant)\to
M_{m,n}(\R_\geqslant)$ be given by
$$(\sigma,(A_1,\ldots,A_n))\longmapsto (A_{\sigma^{-1}(1)},\ldots,A_{\sigma^{-1}(n)}),$$
where $A=(A_1,\ldots,A_n)$ is an arbitrary $m\times n$ matrix
written in a column form. Then $O(n)^\tro$ is given by the orbit of
the identity matrix under the action  $\Sigma_n\times
M_{n,n}(\R_\geqslant)\to M_{n,n}(\R_\geqslant)$.
\end{lmm}
Notice that there is an inclusion, in fact a map of monoids,
$$O(n)^\tro\longrightarrow Gl(\R_\geqslant,n).$$

The above lemma tells us that the action of $O(n)^\tro\times
M_{n,n}(\R_\geqslant)\to M_{n,n}(\R_\geqslant)$ given by the
tropical matrix multiplication, will be only the permutation of the
rows of a given matrix. Let us write $GL(\R_\geqslant,n)$ for the
set of all tropical $n\times n$ invertible matrices. For $A\in
GL(\R_\geqslant,n+k)$ we write
$$A=\left(\begin{array}{ll}
          A_{nn} & B\\
          C      & A_{kk}\end{array}\right)$$
where $A_{nn}$ is the $n\times n$ and $A_{kk}$ is a $k\times k$
block. This allows one to define the action of $O(n)^\tro\times
O(k)^\tro$ on $GL(\R_\geqslant,n+k)$. One then will guess that there
should be a one to one correspondence
$$\frac{Gl(\R_\geqslant,n+k)}{O(n)^\tro\times O(k)^\tro}\longrightarrow G_k(\R^{n+k}_\geqslant).$$
Finally, notice that $O(n)^\tro$ is a monoid under the tropical
matrix multiplication.

\subsubsection{Idempotents}
Suppose $M$ is an arbitrary monoid. An element $a\in M$ is
idempotent if $a^2=a$. We consider to the problem of determining the
idempotent in the monoid $M_{n,n}(\R_\geqslant)$. The result reads
as following.
\begin{lmm}
Let $A$ be an idempotent $n\times n$ matrix entries from
$\R_\geqslant$. Then $A$ satisfies the following conditions
$$\begin{array}{lll}
A_{ii}&\leqslant 1\\
A_{ik}A_{ki}&\leqslant\min(A_{ii},A_{kk}) &\textrm{if }i\neq k.
\end{array}$$
\end{lmm}
The proof is straightforward once we consider the diagonal elements.
The equation $A\odot A=A$ implies that
$$\max(A_{ik}A_{ki}:1\leqslant k\leqslant n)=A_{ii}.$$
For instance, this implies that $A_{ii}^2\leqslant A_{ii}$ which
means that $A_{ii}\leqslant 1$. The other inequality is obtained in
a similar fashion, by comparing the equations for $A_{ii}$ and
$A_{kk}$.

\subsubsection{Stablisation}
We consider the idea of the infinite dimensional orthogonal tropical
matrices. Notice that for any $m,n$ there is a mapping
$M_{m,n}(\R_\geqslant)\to M_{m+1,n+1}(\R_\geqslant)$ given by
$$A\longmapsto\left(\begin{array}{ll}
                    A & 0\\
                    0 & 1\end{array}\right).$$
In particular, this induces a mapping $O(n)^\tro\to O(n+1)^\tro$. We
then define the analogous of the infinite orthogonal group by
$$O^\tro=\colim O(n)^\tro.$$
This object inherits a monoid structure induces from the monoid
structure on the finite dimensional tropical orthogonal monoids. One
then hopes that this will give a characterisation of the infinite
dimensional Grassmannians.

\subsubsection{Comments on Tropical Bundle Theory}
The algebraic topology of fibre bundles with a given topological
space $F$ as the fibre, is understood in terms of the classifying
space of the groups of automorphisms of $F$.\\
By analogy one may consider to fibre bundle theory of surjections
$E\to B$ whose fibres are given by copies of the tropical space
$\R^n_\geqslant$. This then makes it quite reasonable to consider
the classifying space $BO(n)^\tro$ of the tropical orthogonal
monoids $O(n)^\tro$ where these are monoids under the tropical
matrix multiplication. The classifying space functor is defined for
monoid (in fact for topological monoids which carry a weaker
structure). Hence, one may observe that the $\R^n_\geqslant$-bundles
are classified in terms of mapping into $BO(n)^\tro$.\\

The interest in such theory, and the theory of associated
characteristic classes, seems to come from the theory of
singularities. We note that a canonical example for a tropical
bundle will be ``tangent bundle'' of a tropical manifold. This then
motivates one to claim that the classification of these
singularities might be done in terms of the characteristic classes
of the associated tangent bundle.\\ \\
Hadi Zare, The School of Mathematics, Manchester University, Manchester, UK M13 9PL, \textit{email: hzare@maths.manchester.ac.uk}

\end{document}